# Fractal algebraic and linguistic structures generated by the Laurent series for the Gamma function near its negative poles

Andrei Vieru

**Abstract**
We give closed-form expressions for the Laurent series coefficients of the Gamma function near all its strictly negative singularities. These closed-form expressions are clearly self-similar. We briefly describe their algebraic and grammatical budding patterns. As the degree of the coefficient grows to infinity its global structure becomes a fractal.

## 1. Introduction and basic notations

Harmonic numbers are defined as

$$H_n = \sum_{k=1}^{n} \frac{1}{k} = \int_0^1 \frac{1-x^n}{1-x} dx \tag{1}$$

Euler's constant is defined as

$$\gamma = \lim_{n \to \infty} \left[ \sum_{k=1}^{n} \frac{1}{k} - \ln(n) \right] \tag{2}$$

Euler's limit definition of the Gamma function

$$\Gamma(z) = \lim_{n \to \infty} \frac{n!\, n^z}{z(z+1)\ldots(z+n)}$$

works directly in the whole complex domain without need of analytic continuation.

The relation of the Euler constant to the Gamma function is well-known, either via the digamma function):

$$-\gamma = \psi(1) = \frac{\Gamma'(1)}{\Gamma(1)}$$

or directly, as a limit formula (which yields also the constant term in the Laurent expansion of Γ near 0):

$$-\gamma = \lim_{x \to 0} \left[ \Gamma(x) - \frac{1}{x} \right] \qquad (3)$$

Euler's constant is related to the set of Harmonic numbers by its original definition due to Euler himself **(2)**, and is related to each Harmonic number in particular by the formula:

$$\psi(n) = H_{n-1} - \gamma$$

It is also related to the set of fractional Harmonic numbers by the formula:

$$\gamma = \int_0^1 H_x \, dx$$

The Riemann zeta function is defined as $\quad \zeta(s) = \sum_{k=1}^{\infty} \frac{1}{k^s}$

## 2. In search of Harmonic numbers near the Gamma function poles

The neighborhoods of zero have been evoked here above.
Inspecting $\Gamma$ near the strictly negative poles, one can find a limit formula which in its general form reads:

$$\lim_{x \to 0} \left[ \frac{(-1)^{n+1}}{n!x} + \Gamma(-n+x) \right] = \frac{(-1)^{n+1}}{n!} (\gamma - H_n) \qquad (4)$$

(for *n* > 0 and with *n*! = $\Gamma(n+1)$)

or

$$\gamma = (-1)^{n+1} (n!) \lim_{x \to 0} \left[ \frac{(-1)^{n+1}}{n!x} + \Gamma(-n+x) \right] + H_n \qquad (n \geq 0)$$

which generalizes (2), considering that $H_0$ = 0 (by virtue of Euler's integral representation, see (1))

## 3. The singularities of Gamma at strictly negative arguments

For *n* = 0, we have the following limit formula:

$$\lim_{x \to 0} \left\{ \frac{\frac{1}{x} - \Gamma(x)}{-1 + x\Gamma(-x)} \right\} = -\frac{\gamma}{2} \qquad (\clubsuit)$$

When we approach strictly negative points of singularity of the Γ function, Euler's constant escorted by Harmonic numbers arise also in the following way (for *n* ≥ 1):

$$\lim_{x \to 0} \left\{ \frac{1}{x} \left[ 1 + \frac{\frac{1}{x} - \Gamma(-n+x)}{-\frac{1}{x} - \Gamma(-n-x)} \right] \right\} = (-1)^{n+1} \frac{2}{n! + (-1)^{n+1}} (\gamma - H_n)$$

$$(\spadesuit)$$

Here, surprisingly, occurs an element of self-similarity: if we divide the expression under limit in the LHS of (♠) by the last factor of the RHS, namely by $\gamma - H_n$ and then subtract the remainder of the RHS from it, we get, multiplying by the denominator in the RHS and dividing by *x*, exactly the same limit:

$$\lim_{x \to 0} \left\{ \frac{n! + (-1)^{n+1}}{(-1)^n x} \left[ \frac{1}{x(\gamma - H_n)} \left( 1 + \frac{\frac{1}{x} - \Gamma(-n+x)}{-\frac{1}{x} - \Gamma(-n-x)} \right) - (-1)^{n+1} \frac{2}{n! + (-1)^{n+1}} \right] \right\}$$

$$= (-1)^{n+1} \frac{2}{n! + (-1)^{n+1}} (\gamma - H_n)$$

One can write (♠) and (♣) as one single formula:

$$\lim_{x \to \infty} \left[ \frac{\chi_{n0}}{x} + \frac{\frac{1}{x} - \Gamma(-n+x)}{-1 - (-1)^{\delta_{n0}} x\Gamma(-n-x)} \right] = (-1)^{n+1} \left( \frac{2}{n! + \chi_{n0}(-1)^{n+1}} \right)^{\varepsilon_{n0}} (\gamma - H_n)$$

where $\chi_{ij} = \begin{cases} 0 & \text{if } i=j \\ 1 & \text{if } i \neq j \end{cases}$, $\delta_{ij}$ is the Kronecker delta symbol and $\varepsilon_{ij} = \begin{cases} -1 & \text{if } i=j \\ 1 & \text{if } i \neq j \end{cases}$

(in the context $\chi_{n0} = \begin{cases} 0 & \text{if } n=0 \\ 1 & \text{if } n \neq 0 \end{cases}$ is equivalent to the Dirac measure of the set **N⁺**)

The need of Kronecker-like symbols in the unified formula points to a specificity of the rightmost point of singularity (zero).

# 4. The closed form expressions of the Laurent series coefficients in the vicinities of the strictly negative poles of the Gamma function

In (4°), page 297, the authors comment Ramanujan's attempt to compute the coefficients of the Taylor series for $\Gamma(1+x)$ and show to which extent Ramanujan's computations are accurate. The authors propose more accurate values than those found in Ramanujan's second notebook, but they manage to print still one or two inaccurate digits in two of their 'corrected' coefficients.

The Laurent expansion of the Gamma function near 0 is well-known:

$$\Gamma(z) = \frac{1}{z} - \gamma + \left(\frac{\gamma^2 + \zeta(2)}{2}\right) z - \frac{2\gamma^3 + \gamma\pi^2 + 4\zeta(3)}{12} z^2 + \ldots \quad \textbf{(5)}$$

The coefficients in **(5)** are related to one another by the following well-known recurrence formula[1], involving the values of the Riemann zeta function at integer arguments, valid for $n > 2$, after choosing $a_1 = 1$ and $a_2 = -\gamma$:

$$(n-1)a_n = -\gamma a_{n-1} + \zeta(2)a_{n-2} - \zeta(3)a_{n-3} + \ldots + (-1)^n \zeta(n) \quad \textbf{(6)}$$

As we have already seen, the rightmost singularity is somewhat different from all other strictly negative poles of the Gamma function. Laurent series near all other singularities inherit the coefficients found in **(5)** using **(6)**, *but only as a part of their own closed-form.*
The beginning of the Laurent expansion near any pole $-n < 0$ reads:

$$\Gamma(-n+z) = \frac{(-1)^n}{n!z} + \frac{(-1)^{n+1}}{n!}(\gamma - H_n) + \frac{(-1)^{n+2}}{n!}\left(\frac{\gamma^2 + \zeta(2)}{2} - \sum_{k=1}^{n}\frac{\gamma - H_k}{k}\right)z$$

$$+ \frac{(-1)^{n+3}}{n!}\left(\frac{2\gamma^3 + \gamma\pi^2 + 4\zeta(3)}{12} - \sum_{j=1}^{n}\frac{\frac{\gamma^2+\zeta(2)}{2} - \sum_{k=1}^{j}\frac{\gamma-H_k}{k}}{j}\right)z^2 + \ldots \quad \textbf{(7)}$$

One immediately can see that the coefficients in **(7)** are exactly those of **(5),** from which a certain summation is subtracted. In the constant term of **(7)** we see $H_n$, and we know from **(1)** that Harmonic numbers *are already* summations. (We shall leave

---

[1] This relation appears in the *Taylor* expansion of $1/\Gamma$ and of $\Gamma(1+x)$ with $na_n$ in the LHS and with indexes shifted in the RHS (and therefore with changed signs, see (7°), p.66)

out the factorials from this discussion: they might have been written as a common factor for all the infinite sum.)

## 5. On the self-similarity of the Laurent series coefficients expressed in closed form

On the other hand, the constant term $\gamma - H_n$, gives birth, so to say, to an efflorescence of nested summations in the coefficients of **(7)**. Since we decided to drag outside the signs of the summands writing them as $(-1)^{n+k}$, let us now consider the absolute value of the coefficients $a_k$ in **(6)** and of the coefficients $b_k$ in **(7)**.

We have $b_1 = 1 = |a_1|$ and $\quad b_2 = \gamma - H_n = |a_2| - \sum_{k=1}^{n} \frac{|a_1|}{k}$ **(8)**

then

$$b_3 = |a_3| - \sum_{k_2=1}^{n} \frac{|a_2| - \sum_{k_1=1}^{k_2} \frac{|a_1|}{k_1}}{k_2}$$ **(9)**

$$b_4 = |a_4| - \sum_{k_3=1}^{n} \frac{|a_3| - \sum_{k_2=1}^{k_3} \frac{|a_2| - \sum_{k_1=1}^{k_2} \frac{|a_1|}{k_1}}{k_2}}{k_3}$$ **(10)**

$$b_5 = |a_5| - \sum_{k_4=1}^{n} \frac{|a_4| - \sum_{k_3=1}^{k_4} \frac{|a_3| - \sum_{k_2=1}^{k_3} \frac{|a_2| - \sum_{k_1=1}^{k_2} \frac{|a_1|}{k_1}}{k_2}}{k_3}}{k_4}$$ **(11)**

etc.

We have a recurrence relation for $n > 2$ which in fact is very easy to understand but impossible to write using standard symbols. The $k_i$ are not only *indexed* indexes (there are not enough letters in a finite alphabet if one wants to write further these coefficients) but also *nested* indexes since, around the pole $-n$, they eventually all run from 1 to $n$, *but depending on each other subordinately*: what is important to stress on is that the current value of $k_i$ in the denominators is the same as the upper bound of the largest summation in the corresponding numerators.

One may doubt whether these closed-form expressions might be useful for very precise computations, although they work quite well for little *z* values. In contrast, they have a nice and very unusual fractal algebraic structure.

The impossibility to write a recurrence relation of the form $b_i = r(b_{i-1})$ stems from the need to change in the $b_{i-1}$ at each step of the recurrence the upper bound of the largest summation in it.

In exchange one can indeed write the recurrence relation in terms of fully written RHS of **(8), (9), (10), (11)**, etc.

For example, if $A_i$ ($i > 1$) is the string written in TeX corresponding to $b_i$ then to obtain the string corresponding to $b_{i+1}$ one has to carry out in the given order the two following patterns of rewriting rules[1]:

n→{k}_{i}
$A_i$→\left|{a}_{i+1}\right|-\sum_{{k}_{i}=1}^{n}{\frac{$A_i$}{{k}_{i}}}

Since $A_i$ is supposed to be a fully written formula, in both rules "i" and "i+1" are supposed to be written in digits. (This is not necessarily the case of "n", which may designate any pole of Gamma.)

The "axiom" is $A_2$:
\left|{a}_{2}\right|-\sum_{{k}_{1}=1}^{n}{\frac{\left|{a}_{1}\right|}{{k}_{1}}}

wich corresponds to $\quad |a_2| - \sum_{k_1=1}^{n} \frac{|a_1|}{k_1}$

applying
n→{ k }_{ 2 } and then
$A_2$→\left|{a}_{3}\right|-\sum_{{k}_{2}=1}^{n}{\frac{$A_2$}{{k}_{2}}}

one gets
\left|{a}_{3}\right|-\sum_{{k}_{2}=1}^{n}{\frac{\left|{a}_{2}\right|-\sum_{{k}_{1}=1}^{{ k }_{ 2 }}{\frac{\left|{a}_{1}\right|}{{k}_{1}}}}{{k}_{2}}}

---

[1] we say "patterns of rules" according to the old John von Neumann's idea for formal systems: in fact we have here an infinity of rules corresponding to all "i" and all $A_i$ that actually need to be replaced by numbers and, respectively, by formulae written in TeX.

which corresponds to $\quad |a_3| - \sum_{k_2=1}^{n} \dfrac{|a_2| - \sum_{k_1=1}^{k_2} \frac{|a_1|}{k_1}}{k_2}$

An example:
1/(2*(10^-7))-(γ-3/2)/2+(1/2)*((γ^2+zeta(2))/2-(γ-1)/1-(γ-3/2)/2)*(10^-7)-(1/2)*((2*γ^3+γ*π^2+4*zeta(3))/12-(((γ^2+zeta(2))/2-(γ-1)/1)+((γ^2+zeta(2))/2-((γ-1)/1+(γ-3/2)/2))/2))*(10^-7)^2+(1/2)*(((2*γ^4+γ^2*π^2+γ*4*zeta(3))/12+zeta(2)*(γ^2+zeta(2))/2+zeta(3)*γ+zeta(4))/4-(((2*γ^3+γ*π^2+4*zeta(3))/12-(((γ^2+zeta(2))/2-(γ-1)/1)/1))/1+((2*γ^3+γ*π^2+4*zeta(3))/12-((((γ^2+zeta(2))/2-(γ-1)/1)/1)/1+((γ^2+zeta(2))/2-((γ-1)/1+(γ-3/2)/2))/2))/2))*(10^-7)^3=
5000000,4613922612108656733690271352 3678…    **(12)**

This truncated expansion of Γ(-2+10⁻⁷) is correct to the 27th decimal.

The coefficient of (10^-7)^3 corresponds to

$$b_5 = |a_5| - \sum_{k_4=1}^{2} \dfrac{|a_4| - \sum_{k_3=1}^{k_4} \dfrac{|a_3| - \sum_{k_2=1}^{k_3} \dfrac{|a_2| - \sum_{k_1=1}^{k_2} \frac{|a_1|}{k_1}}{k_2}}{k_3}}{k_4}$$

The fractal features of **(12)** are visible with the naked eye, especially in fractions. Parts of the expression like (γ-1)/1)/1)/1 are of course useless — unless, for the sake of the example, we want to write the full string "grammatically".

Andrei Vieru
andreivieru007@gmail.com